\documentclass[twocolumn]{article}
\usepackage{amstext,amssymb,amsbsy,array,amsmath,graphicx}
\usepackage{graphicx}
\usepackage{stmaryrd}

\newcommand{\nsd}{{n_{\text{sd}}}}

\newcommand{\bfv}{{\boldsymbol v}}
\newcommand{\bfp}{{\boldsymbol p}}

\newcommand{\bfsig}{{\boldsymbol\sigma}}
\newcommand{\bfh}{{\boldsymbol h}}

\newcommand{\bfq}{{\boldsymbol q}}

\newcommand{\bfI}{{\boldsymbol 1}}
\newcommand{\bfu}{\boldsymbol u}

\newcommand{\bfn}{{\boldsymbol n}}

\newcommand{\bftau}{{\boldsymbol\tau}}
\newcommand{\bfeps}{{\boldsymbol\varepsilon}}

\newcommand{\bff}{{\boldsymbol f}}

\newcommand{\bfg}{{\boldsymbol g}}

\begin{document}
\title{\bf A mixed method for elasticity with the curl of displacements as a drilling degree of freedom}
%
\author{Peter Hansbo\\ Department of Mechanical Engineering,\\ J\"onk\"oping University,\\
SE-55111 J\"onk\"oping, Sweden}
%

\maketitle

\begin{abstract}
{We present a mixed method for the linearized elasticity equations with independent approximation of the curl of the displacements.
The curl can be seen as a drilling degree of freedom allowing for coupling with rotating objects and the direct application of moments of force.
}
\end{abstract}

%

\section{Introduction}

A drilling degree of freedom refers to a rotational degree of freedom for membrane elements, e.g. for
linearized elasticity. 
Such degrees of freedoms refer to forces rather than displacements and are thus not natural to
incorporate as enhancements of displacement fields. Nevertheless, the early attempts at formulating elements
for planar elasticity with drilling degrees of rotations, these were seen as enhancements of the displacement field, see, e.g., 
Allman \cite{All84,All88} and Bergan and Felippa \cite{BeFe85}. Such elements have to be carefully constructed using
particular qualities of the basic elements they are meant to enhance, they are difficult to generalize to higher order, and are plagued by stability problems.
Since the drilling degree of freedom is a force-type variable, mixed methods would seem more natural and were
considered by Hughes and Brezzi \cite{HuBr89}, where several different methods were proposed. The simplest of these was subsequently studied
from a numerical point of view by Hughes, Masud, and Harari \cite{HuMaHa95}. 

In this note, a different approach to mixed methods for drilling degrees of freedom is taken for the case of isotropic linearized elasticity. Instead of artificially
adding the drilling degree of freedom to the energy functional, as in \cite{HuMaHa95}, which can lead to numerical stiffening
of the discrete problem, we introduce the drilling degree of freedom by splitting the weak form of the elasticity equations 
in a suitable way so that 
the curl of displacements is identified and can be replaced by an independent variable. 

The remainder is organized as follows: in Section 2 we introduce the model problem and the split of the equations ; in Section 3 we remark on some different choices of finite element spaces; in Section 4 we compare and contrast our method with that of \cite{HuMaHa95}; in Section 5
we present some numerical experiments, in particular comparing with \cite{HuMaHa95} and with standard conforming methods.

\section{The Continuous Problem}

Consider a domain $\Omega$ in ${\mathbb{R}}^\nsd$, $\nsd=2$ or $\nsd=3$ with boundary
$\Gamma = \Gamma_\text{D}\cup\Gamma_\text{N}$, $\Gamma_\text{D}\cap\Gamma_\text{N} =\emptyset$, whose outward pointing normal is denoted
$\bfn$. The linear elasticity equations can be written
\begin{equation}\begin{array}{rcl} \label{strongform}
-\nabla\cdot\bfsig (\bfu) & = & {\boldsymbol f}\;\quad\text{in}\;\Omega\\[3mm]
\bfu & = & {\bf 0}\quad\text{on}\;\Gamma_\text{D}\\[3mm]
\bfn\cdot\bfsig(\bfu) & = & \bfg\quad\text{on}\;\Gamma_\text{N}
\end{array}\end{equation}
Here, with $\lambda$ and $\mu$ given material data, the stress tensor $\bfsig$ is
defined by
\begin{equation}
\bfsig(\bfu)  =  2\mu\,\bfeps(\boldsymbol u) +\lambda\,  \nabla\cdot\bfu\, \bfI
\end{equation}
where $\boldsymbol u$ is the displacement field, $\bfI$ is the identity tensor,
\begin{equation*}
{\boldsymbol\varepsilon }(\boldsymbol u)=\frac 12\left( \nabla \otimes %
\boldsymbol u+(\nabla \otimes \boldsymbol u)^{\text{T}}\right)
\end{equation*}
is the strain tensor, and ${%
\boldsymbol f}$  and $ \bfg$ are given data. We have also used the notation
\[
\left(\nabla\cdot\bftau \right)_i= \sum_{j=1}^\nsd\frac{\partial \tau_{ij}}{\partial x_j}
\]
for the divergence of a tensor field $\bftau$.
Introducing the Hilbert space
\[
V = \{\bfv\in [H^1(\Omega)]^\nsd : \bfv = {\bf 0}\;\;\text{on}\;\;\Gamma_\text{D}\},
\]
the weak form of the elasticity equations is to find $\bfu\in V$ such that
\begin{align}\nonumber
\int_{\Omega}2\mu \bfeps(\bfu) : \bfeps(\bfv)\, dx + \int_{\Omega}\lambda \nabla\cdot\bfu\,\nabla\cdot\bfv \, dx \\
= \int_{\Omega}\bff\cdot\bfv\, dx + \int_{\Gamma_\text{N}}\bfg\cdot\bfv \, ds\label{weakform}
\end{align}
for all $\bfv\in V$. Here $\bfsig:\bfeps = \sum_{ij}\sigma_{ij}\varepsilon_{ij}$. 

The idea is to now rearrange the bilinear form in such a way that we can isolate
the term 
\[
\int_{\Omega}\mu \nabla\times \bfu \cdot\nabla\times \bfv \, dx
\]
which corresponds to the rotational part of the form. It is also symmetric, leaving the remainder symmetric in turn.
We shall then introduce a new variable $\bfp = \nabla\times \bfu$ and construct a mixed method. Note that $\bfp$ is scalar in 2D.

We remark that it is difficult to isolate the curl using the strong form of the equations, for instance using the identity $\nabla (\nabla\cdot\bfu) = \nabla\cdot(\nabla\otimes\bfu)^{\rm T} +\nabla\times\nabla\times\bfu$, since the equations must be on the conservation form (\ref{strongform}) in order to give the right weak boundary conditions in the corresponding weak form (\ref{weakform}). Modifications of the equations are thus best performed on the weak form \emph{after}\/ integration by parts.

\subsection{Rearranging the bilinear form}

In two dimensions we can explicitly write
\begin{align*}
\bfeps(\bfu) : \bfeps(\bfv) = {}& \frac{\partial u_1}{\partial x_1}\frac{\partial v_1}{\partial x_1}+\frac{\partial u_2}{\partial x_2}\frac{\partial v_2}{\partial x_2}\\ & +\frac{1}{2}\left(\frac{\partial u_2}{\partial x_1}
+\frac{\partial u_1}{\partial x_2}\right)\left(\frac{\partial v_2}{\partial x_1}+\frac{\partial v_1}{\partial x_2}\right)
\end{align*}
and
\[
 \nabla\times \bfu \cdot \nabla\times \bfv  = \left(\frac{\partial u_2}{\partial x_1}-\frac{\partial u_1}{\partial x_2}\right)\left(\frac{\partial v_2}{\partial x_1}-\frac{\partial v_1}{\partial x_2}\right),
\]
which means that 
\begin{align*}
\bfeps(\bfu) : \bfeps(\bfv) = {}& \frac{1}{2}\nabla\times \bfu \cdot \nabla\times \bfv 
 + \frac{\partial u_1}{\partial x_1}\frac{\partial v_1}{\partial x_1}\\ & +\frac{\partial u_2}{\partial x_2}\frac{\partial v_2}{\partial x_2}+\frac{\partial u_1}{\partial x_2}\frac{\partial v_2}{\partial x_1}+\frac{\partial u_2}{\partial x_1}\frac{\partial v_1}{\partial x_2},
\end{align*}
or
\begin{align}\nonumber
\bfeps(\bfu) : \bfeps(\bfv) = {}&\frac{1}{2}\nabla\times \bfu \cdot \nabla\times \bfv \\ & + \sum_i\frac{\partial u_i}{\partial x_i}\frac{\partial v_i}{\partial x_i}+\sum_{i\neq j}\frac{\partial u_i}{\partial x_j}\frac{\partial v_j}{\partial x_i}.\label{rearranged}
\end{align}
It is easily checked that relation (\ref{rearranged}) holds also in the three dimensional case.

Introducing the variable $\bfp$, we can write the 3D system on mixed form as the problem of finding $(\bfu,\bfp)\in V\times Q$,
where 
\[
Q=\{\bfv:\;\bfv\in [L_2(\Omega)]^3\},
\]
such that
\begin{equation}\label{mixedform1}
a(\bfu,\bfv) +b(\bfp,\bfv) = f(\bfv),\quad\forall \bfv\in V
\end{equation}
and
\begin{equation}\label{mixedform2}
b(\bfq,\bfu) - c(\bfp,\bfq) = g(\bfq),\quad\forall \bfq\in Q
\end{equation}
where
\begin{align*}
a(\bfu,\bfv) := {}& \int_{\Omega}2\mu\left(\sum_i\frac{\partial u_i}{\partial x_i}\frac{\partial v_i}{\partial x_i}+\sum_{i\neq j}\frac{\partial u_i}{\partial x_j}\frac{\partial v_j}{\partial x_i}\right)\, dx \\ & + \int_{\Omega}\lambda \nabla\cdot\bfu\,\nabla\cdot\bfv \, dx,
\end{align*}
\[
b(\bfp,\bfv) := \int_{\Omega}\mu \bfp\cdot\nabla\times\bfv \, dx,
\]
\[
c(\bfp,\bfq) := \int_{\Omega}\,  \bfp\cdot\bfq \, dx,
\]
and by Helmholtz decomposition, cf. \cite{GiRa86}, we split the force field $\bff := \nabla\phi - \nabla\times\bfh$ into a gradient field and a rotational field
(with $\bfh\times\bfn = {\bold 0}$ on $\Gamma_\text{N}$) and let
\[
f(\bfv) := \int_{\Omega}\nabla\phi\cdot\bfv\, dx + \int_{\Gamma_\text{N}}\bfg\cdot\bfv \, ds,
\]
and 
\[
g(\bfq) = \int_{\Omega}\bfh\cdot\bfq\, dx .
\]
The split of $\bff$ is not necessary but allows us to directly apply distributed moments as external loads.

By formally solving the second equation, we get $\bfp =  \mu\nabla\times\bfu - \bfh$ and from the first equation then formally follows
\[
-\nabla\cdot\bfsig (\bfu)  =  \nabla\phi -\nabla\times\bfh = \bff .
\]

In the 2D case the only difference is that $\bfp$ and $\bfh$ are replaced by scalar fields, $Q=\{v:\;v\in L_2(\Omega)\}$.

\section{Finite element approximation}

We replace the continuous spaces by discrete counterparts $Q^h\subset Q$, $V^h\subset V$ and pose the problem of finding
$(\bfu^h,\bfp^h)\in V^h\times Q^h$ such that
\begin{equation}\label{mixedformd1}
a(\bfu^h,\bfv) +b(\bfp^h,\bfv) = f(\bfv),\quad\forall \bfv\in V^h
\end{equation}
and
\begin{equation}\label{mixedformd2}
b(\bfq,\bfu^h) - c(\bfp^h,\bfq) = g(\bfq),\quad\forall \bfq\in Q^h
\end{equation}
The question then arises of what restriction there is on the combination of spaces $V^h$ and $Q^h$ with respect to stability of the discrete problem.

A more common format for mixed methods is the case when $c(\cdot,\cdot)=0$: find
$(\bfu^h,\bfp^h)\in V^h\times Q^h$ such that
\begin{equation}\label{mixedformd3}
a(\bfu^h,\bfv) +b(\bfp^h,\bfv) = f(\bfv),\quad\forall \bfv\in V^h
\end{equation}
and
\begin{equation}\label{mixedformd3}
b(\bfq,\bfu^h)  = g(\bfq),\quad\forall \bfq\in Q^h
\end{equation}
This saddle point problem requires coercivity on the kernel of $b(\cdot,\cdot)$: there exists a constant $\alpha$ such that
\[
a(\bfv,\bfv) \geq \alpha \| \bfv\|_{H^1(\Omega)}\quad \forall \bfv\in K^h,
\]
where
\[
K^h = \{\bfv\in V^h : b(\bfq,\bfv) =0\quad \forall \bfq\in Q^h\},
\]
and an \emph{inf--sup} condition: there exists a constant $\beta$ such that
\[
\inf_{\bfq\in Q^h}\sup_{\bfv\in V^h}\frac{b(\bfq,\bfv)}{\| \bfv\|_{H^1(\Omega)}\|\bfq\|_{L_2(\Omega)}} \geq \beta  .
\]
In such a mixed method, the space $Q^h$ cannot be chosen too \emph{large}\/ in order not to overconstrain the problem, which would lead to stability problems.
The presence of $c(\cdot,\cdot)$ relaxes the \emph{inf--sup}\/ condition, but we still need coercivity on $K^h$.
Note that the bilinear form $a(\cdot,\cdot)$ does not fulfill a Korn-type inequality and is thus not coercive on the whole of $V^h$.
In consequence, we need $Q^h$ to be large enough to get sufficient control of $\nabla\times\bfv$ to regain Korn's inequality on $K^h$.
In our case we will thus have the opposite problem of not letting $Q^h$ be too \emph{small}\/ in order to avoid stability problems.
We note that the choice
\[
Q^h := \nabla_h\times V^h ,
\]
i.e., functions in $Q^h$ are chosen as element-wise curl  ($\nabla_h\times$) of functions in $V^h$, will give a method equivalent to the original elasticity problem, since $\bfp$ can then be eliminated element-wise, and the mixed method
will in that event be stable (for example a piecewise linear approximation in $V^h$ and piecewise constant in $Q^h$). More generally, if
\[
\nabla_h\times V^h \subseteq Q^h
\]
the method will, for the same reason, be stable (though no gain comes from increasing the size of $Q^h$), but other choices will also work. For instance, in our numerical experience, the natural equal-order interpolation is stable. In the numerical examples we will show how some different choices of spaces behave.

\section{An alternative method}

The method of Hughes et al. \cite{HuBr89,HuMaHa95} can, in the setting of using curl as an independent variable, be written in terms
of minimization of a modified energy functional
\begin{align}\nonumber
{\mathcal E}(\bfu,\bfp) ={}& \frac{1}{2}\int_{\Omega}\bfsig(\bfu) :\bfeps(\bfu)\, d\Omega \\ \nonumber
& + \frac{\gamma}{2}\int_{\Omega}\left(\bfp-\nabla\times\bfu\right)\cdot\left(\bfp-\nabla\times\bfu\right)\, d\Omega \\
& -\int_{\Omega}\bff\cdot\bfu\, d\Omega -\int_{\Gamma_\text{N}}\bfg\cdot\bfu\, ds \label{energyfunc}
\end{align}
We note that here the new variable $\bfp$ is here introduced via a penalty--like functional with $\gamma$ as penalty parameter to be chosen. Setting $\gamma = \mu$ as recommended in \cite{HuMaHa95} and minimizing the energy with respect to $(\bfu,\bfp)$, followed by discretization,
leads to the problem of finding  $(\bfu^h,\bfp^h)\in V^h\times Q^h$ such that
\begin{align}\nonumber
\int_{\Omega}\bfsig(\bfu^h):\bfeps(\bfv)\, d\Omega
+\int_{\Omega}\mu\nabla\times\bfu^h\cdot\nabla\times\bfv d\Omega  \\ \label{mixedform1} -b(\bfp^h,\bfv) = f(\bfv),\quad\forall \bfv\in V^h
\end{align}
and
\begin{equation}\label{mixedform2}
b(\bfq,\bfu^h) - c(\bfp^h,\bfq) = g(\bfq),\quad\forall \bfq\in Q^h,
\end{equation}
where we used the same split of the load vector as earlier.
Here we can 
again choose $\nabla_h\times V^h \subseteq Q^h$ and eliminate $\bfp^h$ to obtain equivalence with the elasticity problem discretized in $V^h$.
Other choices have the drawback of leading to an additional stiffening of the discrete problem, as is obvious from the modified energy functional (\ref{energyfunc}).
The method is however coercive on the whole of $V^h$ so that in general we trade stability for accuracy compared with (\ref{mixedformd1}) and (\ref{mixedformd2}). In the following Section, we will compare the performance of the different methods for equal order interpolation on bilinear and linear elements.

\section{Numerical examples}

\subsection{Convergence for different combinations of spaces}

\begin{figure}[h]
\begin{centering}
\includegraphics[height=8.5cm]{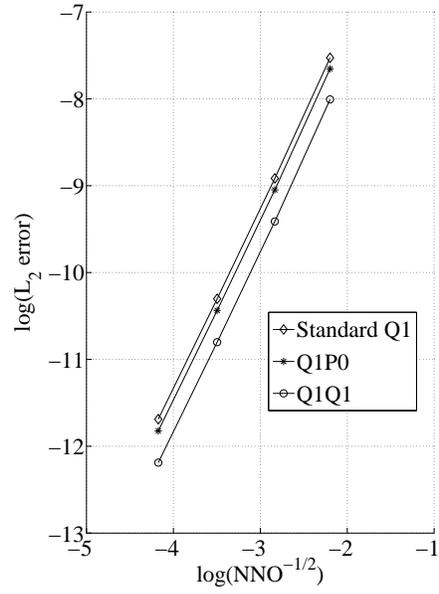}
\end{centering}
\caption{Convergence of bilinear approximations.\label{fig:conv}}
\end{figure}

\begin{figure}[h]
\begin{centering}
\includegraphics[height=8.5cm]{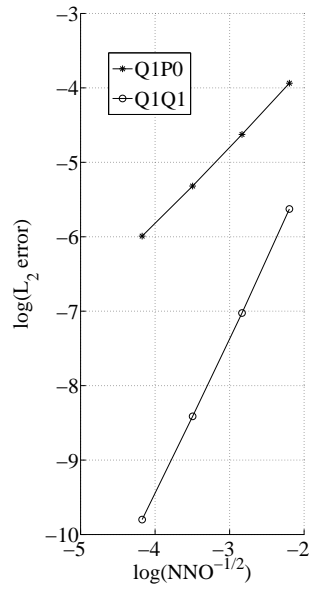}
\end{centering}
\caption{Convergence of the curl.\label{fig:convp}}
\end{figure}

We construct a right-hand side so that the exact solution is 
\[
\bfu=\left[\begin{array}{c}(x_2-1/2)x_2x_1(1-x_1)(1-x_2)\\
-x_2x_1(x_1-1/2)(1-x_1)(1-x_2)\end{array}\right] .
\]
We set $\lambda=\mu=1$ and check the convergence in $L_2$--norm of the Q1Q1 (equal order bilinear approximations) and Q1P0 (bilinear displacements, piecewise constant curl)  and compare with the standard bilinear me\-thod. In Fig. \ref{fig:conv} we show the result which is second order convergence with slightly better error constants for the mixed methods. In Fig. \ref{fig:convp} we show the convergence of the curl variable. Note that, at least on structured meshes, we get second order convergence of the curl for the Q1Q1 element.

\subsection{Comparison with the method of Hughes et al.}

We compare some different methods on an example consisting of a console defined by the domain $0\leq x_1\leq 1$, $0\leq x_2\leq 1$, clamped at $x_1=0$ and with a surface
traction $\bfg = (0,-1)$ at $x_2=1$; no volume load. In plane strain, with a Young's modulus of $E=1$ and Poisson's ratio $\nu=0.3$,
the solution has the approximative ``energy'' 
\[
\| \bfu\|_\sigma^2 := \int_{\Omega}\bfsig(\bfu):\bfeps(\bfu)\, d\Omega \approx 1.903697 , 
\]
as given by Ainsworth et al. \cite{Ains89}.

\begin{figure}[h]
\begin{centering}
\includegraphics[height=8.5cm]{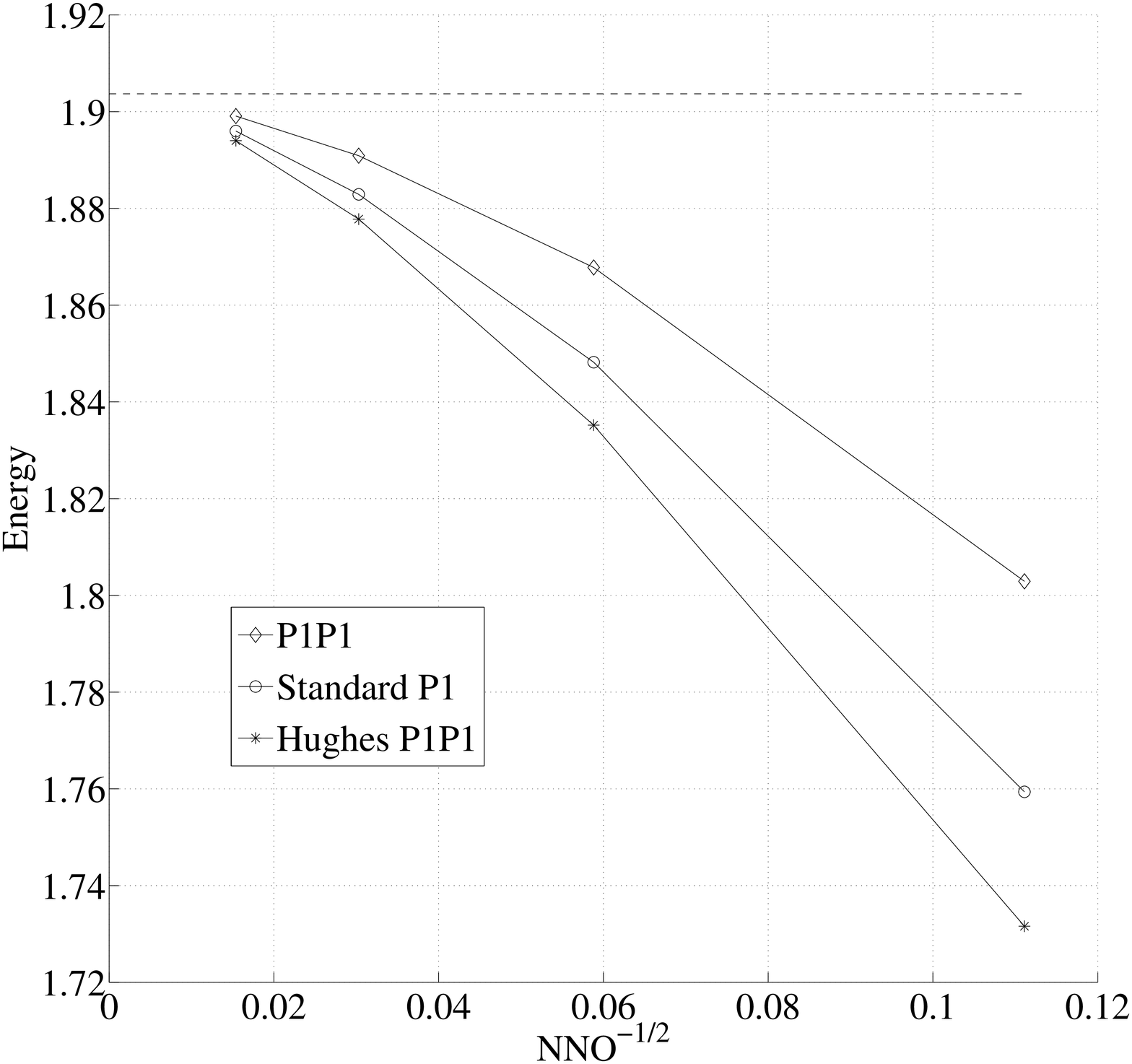}
\end{centering}
\caption{Convergence of different linear methods in ``energy'', $\|\bfu^h\|_\sigma^2$.\label{fig:p1comp}}
\end{figure}

\begin{figure}[h]
\begin{centering}
\includegraphics[height=8.5cm]{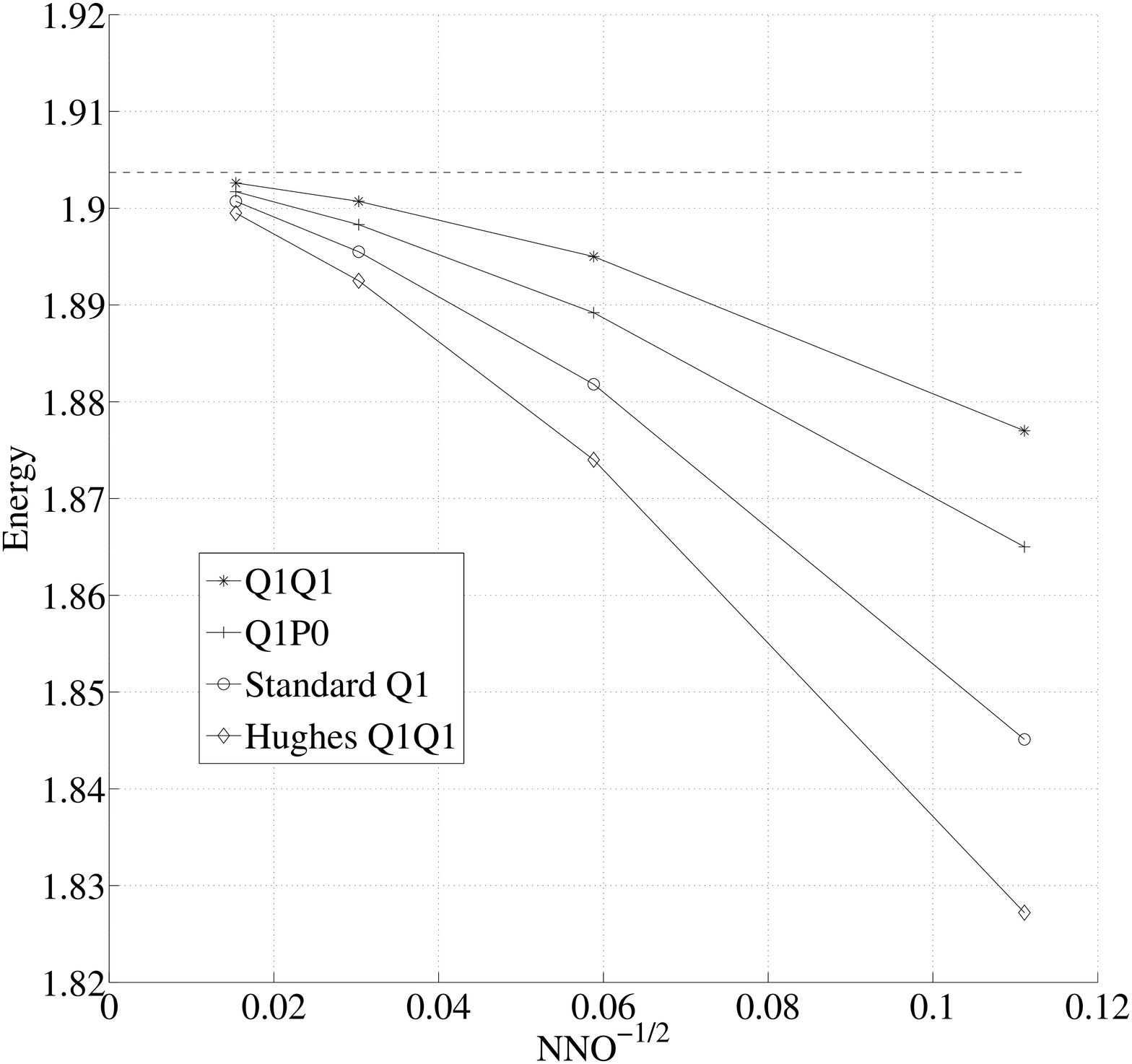}
\end{centering}
\caption{Convergence of different Q1--methods in ``energy'', $\|\bfu^h\|_\sigma^2$.\label{fig:q1comp}}
\end{figure}

In Fig. \ref{fig:p1comp} we compare three different methods: mixed P1P1, standard P1 (constant strain elements) and the method 
(\ref{mixedform1})--(\ref{mixedform2}), which we call ``Hughes method'' in the following. We note that the new method is less stiff than the standard constant strain method, whereas Hughes method adds additional numerical stiffness. The same situation occurs if we
take Q1---elements with equal order interpolation and with piecewise constant approximations for the rotation, see Fig. \ref{fig:q1comp}.
The equal order interpolated method is the least stiff.

\section{Concluding remarks}

We have introduced an approach to drilling degrees of freedom which is close to previously studied methods \cite{HuBr89,HuMaHa95} but which introduces 
no artificial stiffening. The method, which is based on an independent approximation of the curl of displacements, works best with equal order interpolation for displacements and curl.

\bibliographystyle{spmpsci}
\bibliography{drilling}
\newpage

\end{document}